# On a Possible Generalization of Fermat's Last Theorem

Dhananjay P. Mehendale
Sir Parashurambhau College, Tilak Road, Pune–411009
India

## Abstract

This paper proposes a generalized ABC conjecture and assuming its validity settles a generalized version of Fermat's last theorem.

1. **Introduction:** The proof of Fermat's Last Theorem by Andrew Wiles marks the end of a mathematical era. Fermat's last theorem (FLT) was finally settled in the affirmative by Andrew Wiles [1] in the year 1995. The present paper is about stating and proving a possible generalization of FLT. This paper proposes generalized version of ABC conjecture and assuming its validity settles generalized version of Fermat's last theorem.
   The *abc* conjecture, [2], also known as the Oesterlé–Masser conjecture, is a conjecture in number theory, first proposed by Joseph Oesterlé(1988) and David Masser (1985) as an integer analogue of the Mason–Stothers theorem for polynomials. The conjecture is stated in terms of three positive co-prime integers, *a*, *b* and *c* (hence the name), and satisfy *a* + *b* = *c*. If *d* denotes the product of the distinct prime factors of *abc*, the conjecture essentially states that *d* is usually not much smaller than *c*.

2. **A Generalization of FLT:** A possible generalization of Fermat's last theorem can be made as given below.

**Theorem 2.1(Generalized FLT):** For every positive integer $k \geq 2$ there exists a positive integer $g(k) = 2k+2$ such that for every $n \geq g(k)$ the diophantine equation

$$x_1^n + x_2^n + x_3^n + \cdots + x_k^n = z^n$$

has no nontrivial solutions in terms of positive integers $x_i, 1 \leq i \leq k$ and z.

**Remark 2.1:** It is clear to see that the original Fermat's last theorem is a special case of this theorem with $k = 2$ and for this case the theorem is true even for smaller indices. For this case, $g(2) = 6$, but actually the theorem holds even for smaller indices, namely, for the values of *n* = 3, 4, and 5 also.

We now proceed to state generalized ABC conjecture. Assuming the validity of this generalized ABC conjecture we actually see that it is straightforward to prove (as a consequence) the above given theorem 2.1.

**A Generalized ABC Conjecture:** We give below three equivalent versions of the generalized ABC conjecture. Using second generalized version of ABC conjecture, namely, generalized ABC conjecture (II) we show that we can obtain generalized FLT.

**Generalized ABC Conjecture (I):** For every $\varepsilon > 0$, there exist only finitely many (k+1)-tuples $(a_1, a_2, \cdots, a_k, b)$ of positive co-prime integers, with $a_1 + a_2 + \cdots + a_k = b$ such that

$$b > rad(a_1 a_2 \cdots a_k b)^{1+\varepsilon}$$

**Definition 2.1:** Radical of *n*, *rad(n)*, denotes product of all distinct prime factors of *n*.

**Illustrations:** *rad*(16) = 2, *rad*(17) = 17, *rad*(18) = 2.3 = 6.

Let $a = 2$, $b = 3^{10} \cdot 109 = 6,436,341$, $c = 23^5 = 6,436,343$, then
$rad(abc) = 15042$.

**Generalized ABC Conjecture (II):** For every $\varepsilon > 0$, there exists a constant $C(\varepsilon)$ such that for all (k+1)-tuples $(a_1, a_2, \cdots, a_k, b)$ of positive co-prime integers, with $a_1 + a_2 + \cdots + a_k = b$ the inequality

$$b < C(\varepsilon) rad(a_1 a_2 \cdots a_k b)^{1+\varepsilon} \text{ holds.}$$

**Definition 2.2:** Quality, $q(a_1, a_2, \cdots, a_k, b)$, of the (k+1)-tuple is given as follows,

$$q(a_1, a_2, \cdots, a_k, b) = \frac{\log(b)}{\log(rad(a_1 a_2 \cdots a_k b))}.$$

**Generalized ABC Conjecture (III):** For every $\varepsilon > 0$, there exist only finitely many (k+1)-tuples $(a_1, a_2, \cdots, a_k, b)$ of positive co-prime integers, with $a_1 + a_2 + \cdots + a_k = b$ such that
$$q(a_1, a_2, \cdots, a_k, b) > 1 + \varepsilon.$$

**Proof of theorem 2.1:** If generalized ABC conjecture (II) is correct when $C(\varepsilon) = 1$ and $\varepsilon = 1$ and when co-prime natural positive numbers exist which are $a_1, a_2, \cdots, a_k, b$ and they satisfy $a_1 + a_2 + \cdots + a_k = b$ then the following inequality holds:
$$b < rad(a_1 a_2 \cdots a_k b)^2$$

We assume the co-prime natural numbers $x_1^n, x_2^n, \cdots, x_k^n, z^n$ which satisfy generalized FLT: i.e. $x_1^n + x_2^n + \cdots + x_k^n = z^n$.

Now, replacing $a_1$ to $x_1^n$, $a_2$ to $x_2^n$, ......, $a_k$ to $x_k^n$, and $b$ to $z^n$ we have (from above inequality) following outcome from equation $x_1^n + x_2^n + \cdots + x_k^n = z^n$ which represents equation for generalized FLT:
we have
$$z^n < rad(x_1^n x_2^n \cdots x_k^n z^n)^2 = rad(x_1 x_2 \cdots x_n z)^2 \leq (x_1 x_2 \cdots x_k z)^2 < (z^{k+1})^2 = z^{2k+2}$$

Therefore, to have a set of co-prime natural numbers $x_1^n, x_2^n, \cdots, x_k^n, z^n$ satisfying $x_1^n + x_2^n + \cdots + x_k^n = z^n$ the index $n$ must satisfy the inequality given below, namely,

$$n < 2k + 2.$$ □

## Acknowledgements

I am thankful to Prof. M. R. Modak and Mr. Vishwas Chiplunkar for useful discussions.

## References

1. Wiles, A. "Modular Elliptic-Curves and Fermat's Last Theorem." *Ann. Math.* **141**, 443-551, 1995.

2. *abc* conjecture: From Wikipedia, the free encyclopedia.